\documentclass{Z}
\DeclareGraphicsExtensions{.pdf,.eps}

\author{Achim Zeileis\\Wirtschaftsuniversit\"at Wien \And
        Gabor Grothendieck}
\Plainauthor{Achim Zeileis, Gabor Grothendieck}

\title{\pkg{zoo}: An \proglang{S3} Class and Methods for
  Indexed Totally Ordered Observations}
\Plaintitle{zoo: An S3 Class and Methods for
  Indexed Totally Ordered Observations}

\Keywords{totally ordered observations, irregular time series,
  regular time series, \proglang{S3}, \proglang{R}}
\Plainkeywords{totally ordered observations, irregular time series,
  regular time series, S3, R}

\Abstract{
  A previous version to this introduction to the \proglang{R} package \pkg{zoo}
  has been published as \cite{zoo:Zeileis+Grothendieck:2005} in the
  \emph{Journal of Statistical Software}.

  \pkg{zoo} is an \proglang{R} package providing an \proglang{S3}
  class with methods for indexed totally ordered observations, such as
  discrete irregular time series. Its key design goals are independence of a
  particular index/time/date class and consistency with base
  \proglang{R} and the \code{"ts"} class for
  regular time series. This paper describes how these are achieved
  within \pkg{zoo} and provides several illustrations 
  of the available methods for \code{"zoo"} objects which include
  plotting, merging and binding, several mathematical operations,
  extracting and replacing data and index, coercion and \code{NA}
  handling. A subclass \code{"zooreg"} embeds regular time series
  into the \code{"zoo"} framework and thus bridges the gap between
  regular and irregular time series classes in \proglang{R}.
}

\begin{document}


\section{Introduction} \label{sec:intro}

The \proglang{R} system for statistical computing
\citep[\url{http://www.R-project.org/}]{zoo:R:2005}
ships with a a class for regularly spaced time series,
\code{"ts"} in package \pkg{stats}, but has no native class for
irregularly spaced time series. With the increased interest in
computational finance with \proglang{R} over the last years
several implementations of classes for irregular time series 
emerged which are aimed particularly at finance applications.
These include the \proglang{S3} classes \code{"timeSeries"}
in package \pkg{fCalendar} from the \pkg{Rmetrics} bundle \citep{zoo:fCalendar:2004}
and \code{"irts"} in package \pkg{tseries} \citep{zoo:tseries:2005}
and the \proglang{S4} class \code{"its"} in package \pkg{its} \citep{zoo:its:2004}.
With these packages available, why would anybody want yet another 
package providing infrastructure for irregular time series?
The above mentioned implementations have in common that they are restricted to a particular
class for the time scale: the former implementation comes with its own time class
\code{"timeDate"} built on top of the \code{"POSIXt"} classes
available in base \proglang{R} whereas the latter two use \code{"POSIXct"} directly.
And this was the starting point for the \pkg{zoo} project: the first author
of the present paper needed
more general support for ordered observations, independent of a particular
index class, for the package \pkg{strucchange}
\citep{zoo:Zeileis+Leisch+Hornik:2002}. Hence, the package was called
\pkg{zoo} which stands for \underline{Z}'s \underline{o}rdered \underline{o}bservations.
Since the first release, a major part of the additions to \pkg{zoo}
were provided by the second author of this paper, so that the name
of the package does not really reflect the authorship anymore.
Nevertheless, independence of a particular index class remained
the most important design goal. While the package evolved to its current
status, a second key design goal became more and more clear: to provide
methods to standard generic functions for the \code{"zoo"} class that 
are similar to those for the \code{"ts"} class (and base \proglang{R} in
general) such that the usage of \pkg{zoo} is very intuitive because
few additional commands have to be learned. 
This paper describes how these design goals are implemented in \pkg{zoo}.
The resulting package provides the \code{"zoo"} class which offers an
extensive (and still growing) set of standard and new methods for working
with indexed observations and `talks' to the classes \code{"ts"}, \code{"its"},
\code{"irts"} and \code{"timeSeries"}. It also bridges the gap
between regular and irregular time series by providing coercion with (virtually)
no loss of information between \code{"ts"} and \code{"zoo"}.
With these tools \pkg{zoo} provides the basic infrastructure for
working with indexed totally ordered observations and the package can be either employed by
users directly or can be a basic ingredient on top of which other more specialized
applications can be built.

The remainder of the paper is organized as follows:
Section~\ref{sec:zoo-class} explains how \code{"zoo"} objects are created
and illustrates how the corresponding methods for plotting, merging and
binding, several mathematical operations, extracting and replacing data
and index, coercion and \code{NA} handling can be used. Section~\ref{sec:combining}
outlines how other packages can build on this basic infrastructure.
Section~\ref{sec:summary} gives a few summarizing remarks and an outlook
on future developments. Finally, an appendix provides a reference card that
gives an overview of the functionality contained in \pkg{zoo}.

\section[The class "zoo" and its methods]{The class \code{"zoo"} and its methods}
\label{sec:zoo-class}

This section describes how \code{"zoo"} series can be created and subsequently
manipulated, visualized, combined or coerced to other classes. In Section~\ref{sec:zoo},
the general class \code{"zoo"} for totally ordered series is described. Subsequently,
in Section~\ref{sec:zooreg}, the subclass \code{"zooreg"} for
regular \code{"zoo"} series, i.e., series which have an index with a specified
frequency, is discussed. The methods illustrated in the remainder of the
section are mostly the same for both \code{"zoo"} and \code{"zooreg"} objects
and hence do not have to be discussed separately. The few differences in merging and
binding are briefly highlighted in Section~\ref{sec:merge}.

\subsection[Creation of "zoo" objects]{Creation of \code{"zoo"} objects}
\label{sec:zoo}

The simple idea for the creation of \code{"zoo"} objects is to have
some vector or matrix of observations \code{x} which are totally ordered
by some index vector. In time series applications, this index is a measure of
time but every other numeric, character or even more abstract vector that
provides a total ordering of the observations is also suitable. Objects
of class \code{"zoo"} are created by the function
\begin{Scode}
zoo(x, order.by)
\end{Scode}
where \code{x} is the vector or matrix of observations\footnote{In principle,
more general objects can be indexed, but currently \pkg{zoo} does not support this.
Development plans are that \pkg{zoo} should eventually support indexed factors,
data frames and lists.} and \code{order.by}
is the index by which the observations should be ordered. It has to be
of the same length as \code{NROW(x)}, i.e., either the same length as \code{x}
for vectors or the same number of rows for matrices.\footnote{The only case
where this restriction is not imposed is for zero-length vectors, i.e., vectors
that only have an index but no data.} The \code{"zoo"} object
created is essentially the vector/matrix as before but has an additional
\code{"index"} attribute in which the index is stored.\footnote{There is some
limited support for indexed factors available in which case the \code{"zoo"}
object also has an attribute \code{"oclass"} with the original class
of \code{x}. This feature is still under development and might change in future
versions.} Both the observations in the vector/matrix \code{x}
and the index \code{order.by} can, in principle, be of arbitrary classes. However, most of the
following methods (plotting, aggregating, mathematical operations) for \code{"zoo"}
objects are typically only useful for numeric observations \code{x}. Special
effort in the design was put into independence from a particular class for
the index vector. In \pkg{zoo}, it is assumed that combination \code{c()},
querying the \code{length()}, value matching \code{MATCH()}, subsetting \code{[,},
and, of course, ordering \code{ORDER()} work when applied to the index. 
In addition, an \code{as.character()} method might improve printed output\footnote{If
an \code{as.character()} method is already defined, but gives not the desired
output for printing, then an \code{index2char()} method can be defined. This is a
generic convenience function used for creating character representations of the
index vector and it defaults to using \code{as.character()}.}
and \code{as.numeric()} could be used for computing distances between indexes, e.g.,
in interpolation. Both methods are not necessary for working with \code{"zoo"} 
objects but could be used if available.
All these methods are available, e.g., for standard numeric and character vectors and for
vectors of classes \code{"Date"}, \code{"POSIXct"} or \code{"times"}
from package \pkg{chron}, but not for the class \code{"dateTime"} in \pkg{fCalendar}.
In the last case, the solution is to provide methods for the above mentioned
functions so that indexing \code{"zoo"} objects with \code{"dateTime"} vectors works
(see Section~\ref{sec:fCalendar} for an example).
To achieve this  independence of the index class, new generic functions for
ordering (\code{ORDER()}) and value matching (\code{MATCH()}) are introduced
as the corresponding base functions \code{order()} and \code{match()} are 
non-generic. The default methods simply call the corresponding base functions, i.e.,
no new method needs to be introduced for a particular index class if the 
non-generic functions \code{order()} and \code{match()} work for this class.

To illustrate the usage of \code{zoo()}, we first load the package and set the
random seed to make the examples in this paper exactly reproducible.

\begin{Schunk}
\begin{Sinput}
R> library(zoo)
R> set.seed(1071)
\end{Sinput}
\end{Schunk}

Then, we create two vectors \code{z1} and \code{z2} with \code{"POSIXct"} 
indexes, one with random observations
\begin{Schunk}
\begin{Sinput}
R> z1.index <- ISOdatetime(2004, rep(1:2, 5), sample(28, 10), 0, 
+     0, 0)
R> z1.data <- rnorm(10)
R> z1 <- zoo(z1.data, z1.index)
\end{Sinput}
\end{Schunk}
and one with a sine wave
\begin{Schunk}
\begin{Sinput}
R> z2.index <- as.POSIXct(paste(2004, rep(1:2, 5), sample(1:28, 
+     10), sep = "-"))
R> z2.data <- sin(2 * 1:10/pi)
R> z2 <- zoo(z2.data, z2.index)
\end{Sinput}
\end{Schunk}
Furthermore, we create a matrix \code{Z} with random observations and a \code{"Date"}
index
\begin{Schunk}
\begin{Sinput}
R> Z.index <- as.Date(sample(12450:12500, 10))
R> Z.data <- matrix(rnorm(30), ncol = 3)
R> colnames(Z.data) <- c("Aa", "Bb", "Cc")
R> Z <- zoo(Z.data, Z.index)
\end{Sinput}
\end{Schunk}
In the examples above, the generation of indexes looks a bit awkward
due to the fact the indexes need to be randomly generated (and there 
are no special functions for random indexes because these are rarely 
needed in practice). In ``real world'' applications, the indexes
are typically part of the raw data set read into \proglang{R} so the
code would be even simpler. See Section~\ref{sec:combining}
for such examples.\footnote{Note, that in the code above a new \code{as.Date}
method, provided in \pkg{zoo}, is used to convert days 
since 1970-01-01 to class \code{"Date"}. See the respective help page 
for more details.}

Methods to several standard generic functions are available for
\code{"zoo"} objects, such as \code{print}, \code{summary}, \code{str}, \code{head},
\code{tail} and \code{[} (subsetting), a few of which are illustrated in
the following.

There are three printing code styles for \code{"zoo"} objects: vectors are by default
printed in \code{"horizontal"} style
\begin{Schunk}
\begin{Sinput}
R> z1
\end{Sinput}
\begin{Soutput}
 2004-01-05  2004-01-14  2004-01-19  2004-01-25  2004-01-27  2004-02-07 
 0.74675994  0.02107873 -0.29823529  0.68625772  1.94078850  1.27384445 
 2004-02-12  2004-02-16  2004-02-20  2004-02-24 
 0.22170438 -2.07607585 -1.78439244 -0.19533304 
\end{Soutput}
\begin{Sinput}
R> z1[3:7]
\end{Sinput}
\begin{Soutput}
2004-01-19 2004-01-25 2004-01-27 2004-02-07 2004-02-12 
-0.2982353  0.6862577  1.9407885  1.2738445  0.2217044 
\end{Soutput}
\end{Schunk}
and matrices in \code{"vertical"} style
\begin{Schunk}
\begin{Sinput}
R> Z
\end{Sinput}
\begin{Soutput}
           Aa          Bb          Cc         
2004-02-02  1.25543390  0.68157316 -0.63292049
2004-02-08 -1.49458326  1.32341223 -1.49442269
2004-02-09 -1.87462247 -0.87329289  0.62733971
2004-02-21 -0.14538608  0.45234903 -0.14597401
2004-02-22  0.22542418  0.53838938  0.23136133
2004-02-29  1.20695518  0.31814222 -0.01129202
2004-03-05 -1.20861025  1.42379785 -0.81614483
2004-03-10 -0.11039563  1.34774254  0.95522468
2004-03-14  0.84202385 -2.73842019  0.23150695
2004-03-20 -0.19019104  0.12308872 -1.51862157
\end{Soutput}
\begin{Sinput}
R> Z[1:3, 2:3]
\end{Sinput}
\begin{Soutput}
           Bb         Cc        
2004-02-02  0.6815732 -0.6329205
2004-02-08  1.3234122 -1.4944227
2004-02-09 -0.8732929  0.6273397
\end{Soutput}
\end{Schunk}
Additionally, there is a \code{"plain"} style which simply first prints the data 
and then the index.

Above, we have illustrated that \code{"zoo"} series can be indexed like vectors
or matrices respectively, i.e., with integers correponding to their observation
number (and column number). But for indexed observations, one would obviously also
like to be able to index with the index class. This is also available in \code{[}
which only uses vector/matrix-type subsetting if its first argument is of class
\code{"numeric"}, \code{"integer"} or \code{"logical"}.

\begin{Schunk}
\begin{Sinput}
R> z1[ISOdatetime(2004, 1, c(14, 25), 0, 0, 0)]
\end{Sinput}
\begin{Soutput}
2004-01-14 2004-01-25 
0.02107873 0.68625772 
\end{Soutput}
\end{Schunk}

If the index class happens to be \code{"numeric"}, the index has to be either insulated in \code{I()}
like \code{z[I(i)]} or the  \code{window()} method can be used (see Section~\ref{sec:window}).

Summaries and most other methods for \code{"zoo"} objects are carried out
column wise, reflecting the rectangular structure. In addition,
a summary of the index is provided.

\begin{Schunk}
\begin{Sinput}
R> summary(z1)
\end{Sinput}
\begin{Soutput}
     Index                           z1          
 Min.   :2004-01-05 00:00:00   Min.   :-2.07608  
 1st Qu.:2004-01-20 12:00:00   1st Qu.:-0.27251  
 Median :2004-02-01 12:00:00   Median : 0.12139  
 Mean   :2004-02-01 09:36:00   Mean   : 0.05364  
 3rd Qu.:2004-02-15 00:00:00   3rd Qu.: 0.73163  
 Max.   :2004-02-24 00:00:00   Max.   : 1.94079  
\end{Soutput}
\begin{Sinput}
R> summary(Z)
\end{Sinput}
\begin{Soutput}
     Index                  Aa                Bb                Cc          
 Min.   :2004-02-02   Min.   :-1.8746   Min.   :-2.7384   Min.   :-1.51862  
 1st Qu.:2004-02-12   1st Qu.:-0.9540   1st Qu.: 0.1719   1st Qu.:-0.77034  
 Median :2004-02-25   Median :-0.1279   Median : 0.4954   Median :-0.07863  
 Mean   :2004-02-25   Mean   :-0.1494   Mean   : 0.2597   Mean   :-0.25739  
 3rd Qu.:2004-03-08   3rd Qu.: 0.6879   3rd Qu.: 1.1630   3rd Qu.: 0.23147  
 Max.   :2004-03-20   Max.   : 1.2554   Max.   : 1.4238   Max.   : 0.95522  
\end{Soutput}
\end{Schunk}

\subsection[Creation of "zooreg" objects]{Creation of \code{"zooreg"} objects}
\label{sec:zooreg}

Strictly regular series are such series observations where the distance between
the indexes of every two adjacent observations is the same. Such series can 
also be described by their frequency, i.e., the reciprocal value of the distance
between two observations. As \code{"zoo"} can be used to store series with arbitrary
type of index, it can, of course, also be used to store series with regular indexes.
So why should this case be given special attention, in particular as there is already
the \code{"ts"} class devoted entirely to regular series? There are two reasons: First,
to be able to convert back and forth between \code{"ts"} and \code{"zoo"}, the frequency
of a certain series needs to be stored on the \code{"zoo"} side. Second, \code{"ts"} is 
limited to strictly regular series and the regularity is lost if some internal observations
are omitted. Series that can be created by omitting some internal observations from strictly
regular series will in the following be refered to as being (weakly) regular.
Therefore, a class that bridges the gap between irregular and strictly regular series
is needed and \code{"zooreg"} fills this gap. Objects of class \code{"zooreg"} inherit
from class \code{"zoo"} but have an additional attribute \code{"frequency"} in which 
the frequency of the series is stored. Therefore, they can be employed to represent
both strictly and weakly regular series.

To create a \code{"zooreg"} object, either the command \code{zoo()} can be used
or the command \code{zooreg()}.

\begin{Scode}
zoo(x, order.by, frequency)
zooreg(data, start, end, frequency, deltat, ts.eps, order.by)
\end{Scode} 

If \code{zoo()} is called as in the previous section but with an additional
\code{frequency} argument, it is checked whether \code{frequency} complies
with the index \code{order.by}: if it does an object of class \code{"zooreg"}
inheriting from \code{"zoo"} is returned. The command \code{zooreg()} takes mostly
the same arguments as \code{ts()}.\footnote{Only if \code{order.by}
is specified in the \code{zooreg()} call, then \code{zoo(x, order.by, frequency)}
is called.} 
In both cases, the index class is more restricted than in the plain \code{"zoo"}
case. The index must be of a class which can be coerced to \code{"numeric"} 
(for checking its regularity) and when converted to numeric 
the index must be expressable as multiples of 1/frequency. 
Furthermore, adding/substracting
a numeric to/from an observation of the index class, should return the correct value
of the index class again, i.e., group generic functions \code{Ops} should be defined.\footnote{An
application of non-numeric indexes for regular series are the classes \code{"yearmon"}
and \code{"yearqtr"} which are designed for monthly and quarterly series respectively
and are discussed in Section~\ref{sec:yearmon}.}

The following calls yield equivalent series

\begin{Schunk}
\begin{Sinput}
R> zr1 <- zooreg(sin(1:9), start = 2000, frequency = 4)
R> zr2 <- zoo(sin(1:9), seq(2000, 2002, by = 1/4), 4)
R> zr1
\end{Sinput}
\begin{Soutput}
   2000(1)    2000(2)    2000(3)    2000(4)    2001(1)    2001(2)    2001(3) 
 0.8414710  0.9092974  0.1411200 -0.7568025 -0.9589243 -0.2794155  0.6569866 
   2001(4)    2002(1) 
 0.9893582  0.4121185 
\end{Soutput}
\begin{Sinput}
R> zr2
\end{Sinput}
\begin{Soutput}
   2000(1)    2000(2)    2000(3)    2000(4)    2001(1)    2001(2)    2001(3) 
 0.8414710  0.9092974  0.1411200 -0.7568025 -0.9589243 -0.2794155  0.6569866 
   2001(4)    2002(1) 
 0.9893582  0.4121185 
\end{Soutput}
\end{Schunk}

to which methods to standard generic functions for regular series can be
applied, such as \code{frequency}, \code{deltat}, \code{cycle}.

As stated above, the advantage of \code{"zooreg"} series is that they remain
regular even if an internal observation is dropped:

\begin{Schunk}
\begin{Sinput}
R> zr1 <- zr1[-c(3, 5)]
R> zr1
\end{Sinput}
\begin{Soutput}
   2000(1)    2000(2)    2000(4)    2001(2)    2001(3)    2001(4)    2002(1) 
 0.8414710  0.9092974 -0.7568025 -0.2794155  0.6569866  0.9893582  0.4121185 
\end{Soutput}
\begin{Sinput}
R> class(zr1)
\end{Sinput}
\begin{Soutput}
[1] "zooreg" "zoo"   
\end{Soutput}
\begin{Sinput}
R> frequency(zr1)
\end{Sinput}
\begin{Soutput}
[1] 4
\end{Soutput}
\end{Schunk}

This facilitates \code{NA} handling significantly compared to \code{"ts"} and makes
\code{"zooreg"} a much more attractive data type, e.g., for time series regression.

\code{zooreg()} can also deal with non-numeric indexes provided that adding \code{"numeric"}
observations to the index class preserves the class and does not coerce to \code{"numeric"}.

\begin{Schunk}
\begin{Sinput}
R> zooreg(1:5, start = as.Date("2005-01-01"))
\end{Sinput}
\begin{Soutput}
2005-01-01 2005-01-02 2005-01-03 2005-01-04 2005-01-05 
         1          2          3          4          5 
\end{Soutput}
\end{Schunk}

To check whether a certain series is (strictly) regular, the new generic function
\code{is.regular(x, strict = FALSE)} can be used:

\begin{Schunk}
\begin{Sinput}
R> is.regular(zr1)
\end{Sinput}
\begin{Soutput}
[1] TRUE
\end{Soutput}
\begin{Sinput}
R> is.regular(zr1, strict = TRUE)
\end{Sinput}
\begin{Soutput}
[1] FALSE
\end{Soutput}
\end{Schunk}

This function (and also the \code{frequency}, \code{deltat} and \code{cycle}) also 
work for \code{"zoo"} objects if the regularity can still be inferred from the data:

\begin{Schunk}
\begin{Sinput}
R> zr1 <- as.zoo(zr1)
R> zr1
\end{Sinput}
\begin{Soutput}
      2000    2000.25    2000.75    2001.25     2001.5    2001.75       2002 
 0.8414710  0.9092974 -0.7568025 -0.2794155  0.6569866  0.9893582  0.4121185 
\end{Soutput}
\begin{Sinput}
R> class(zr1)
\end{Sinput}
\begin{Soutput}
[1] "zoo"
\end{Soutput}
\begin{Sinput}
R> is.regular(zr1)
\end{Sinput}
\begin{Soutput}
[1] TRUE
\end{Soutput}
\begin{Sinput}
R> frequency(zr1)
\end{Sinput}
\begin{Soutput}
[1] 4
\end{Soutput}
\end{Schunk}

Of course, inferring the underlying regularity is not always reliable and it is safer
to store a regular series as a \code{"zooreg"} object if it is intended to be a regular series.

If a weakly regular series is coerced to \code{"ts"} the missing observations are filled
with \code{NA}s (see also Section~\ref{sec:NA}).
For strictly regular series with numeric index, the class can be switched
between \code{"zoo"} and \code{"ts"} without loss of information.

\begin{Schunk}
\begin{Sinput}
R> as.ts(zr1)
\end{Sinput}
\begin{Soutput}
           Qtr1       Qtr2       Qtr3       Qtr4
2000  0.8414710  0.9092974         NA -0.7568025
2001         NA -0.2794155  0.6569866  0.9893582
2002  0.4121185                                 
\end{Soutput}
\begin{Sinput}
R> identical(zr2, as.zoo(as.ts(zr2)))
\end{Sinput}
\begin{Soutput}
[1] TRUE
\end{Soutput}
\end{Schunk}

This enables direct use of functions such as \code{acf}, \code{arima}, \code{stl} etc. on \code{"zooreg"}
objects as these methods coerce to \code{"ts"} first. 
The result only has to be coerced back to \code{"zoo"}, if appropriate.

\subsection{Plotting}
\label{sec:plot}

The \code{plot} method for \code{"zoo"} objects, in particular for
multivariate \code{"zoo"} series, is based on the corresponding
method for (multivariate) regular time series. It relies on \code{plot}
and \code{lines} methods being available for the index class which can
plot the index against the observations.

By default the \code{plot} method creates a panel for each series
\begin{Schunk}
\begin{Sinput}
R> plot(Z)
\end{Sinput}
\end{Schunk}
but can also display all series in a single panel
\begin{Schunk}
\begin{Sinput}
R> plot(Z, plot.type = "single", col = 2:4)
\end{Sinput}
\end{Schunk}

\begin{figure}[b!]
\begin{center}
\includegraphics{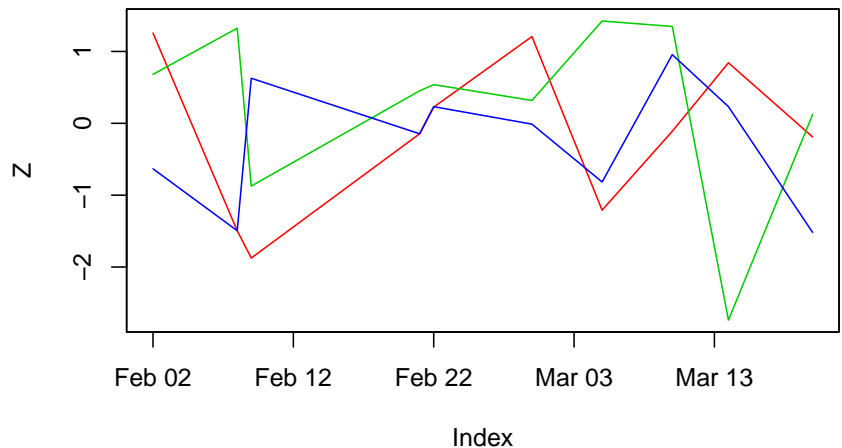}
\caption{\label{fig:plot2} Example of a single panel plot}
\end{center}
\end{figure}

\begin{figure}[p]
\begin{center}
\includegraphics{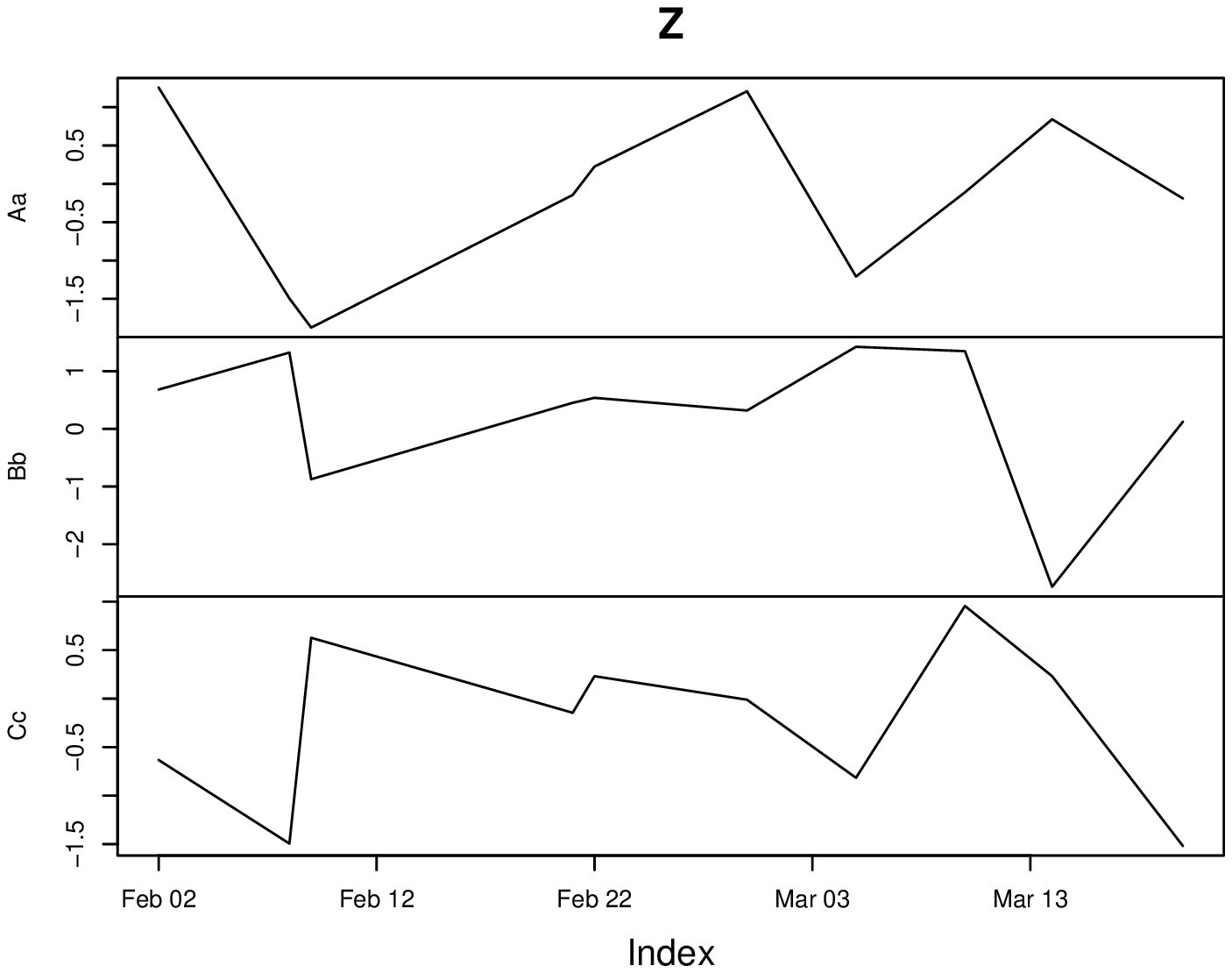}
\includegraphics{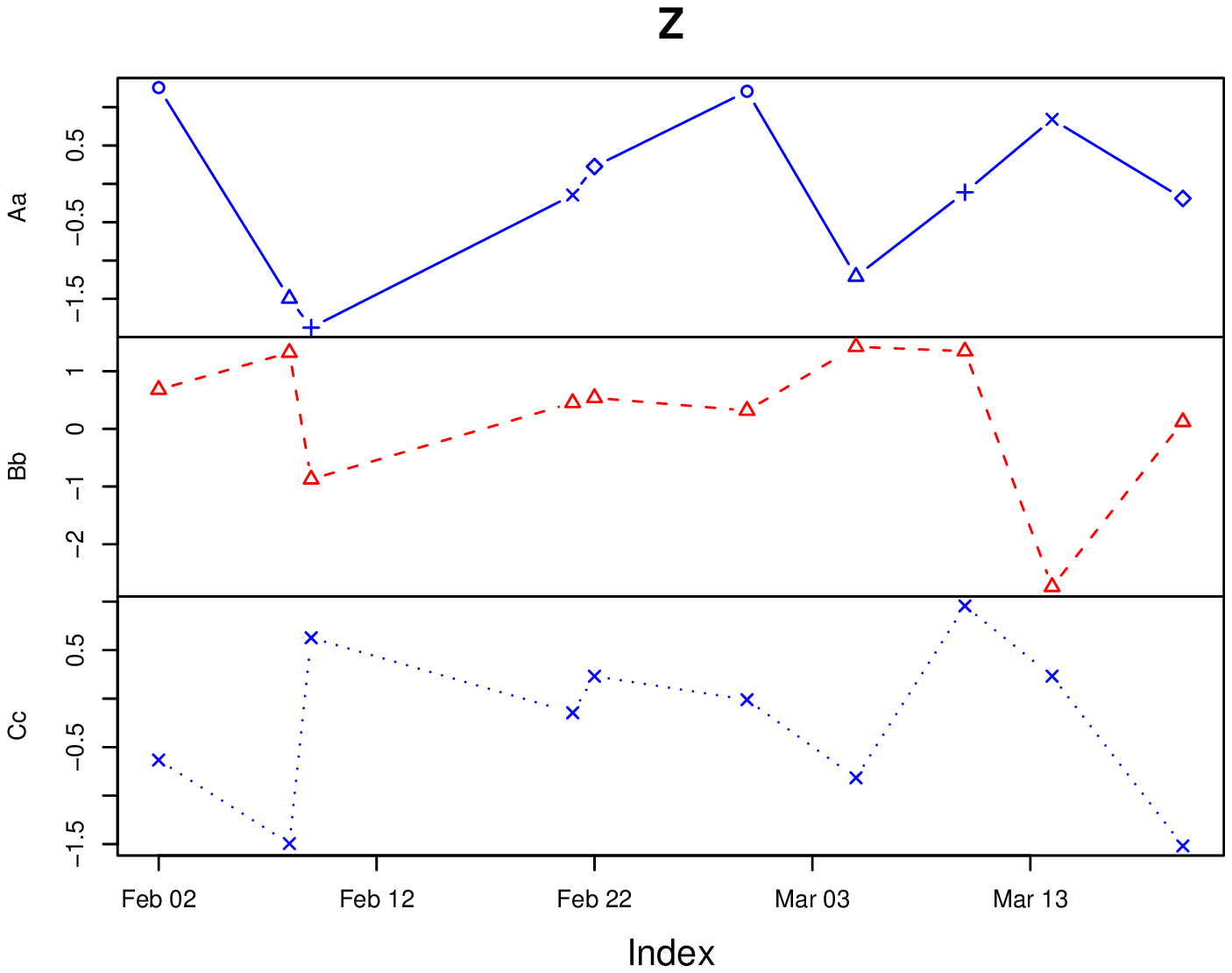}
\caption{\label{fig:plot13} Examples of multiple panel plots}
\end{center}
\end{figure}

In both cases additional graphical parameters like color \code{col},
plotting character \code{pch} and line type \code{lty} can be
expanded to the number of series. But the \code{plot} method for
\code{"zoo"} objects offers some more flexibility in specification
of graphical parameters as in
\begin{Schunk}
\begin{Sinput}
R> plot(Z, type = "b", lty = 1:3, pch = list(Aa = 1:5, Bb = 2, Cc = 4), 
+     col = list(Bb = 2, 4))
\end{Sinput}
\end{Schunk}
The argument \code{lty} behaves as before and sets every series in another
line type. The \code{pch} argument is a named list that assigns to each series
a different vector of plotting characters each of which is expanded to the 
number of observations. Such a list does not necessarily have to include the names of all
series, but can also specify a subset. For the remaining series the default parameter
is then used which can again be changed: e.g., in the above example the \code{col} argument
is set to display the series \code{"Bb"} in red and all remaining series in blue.
The results of the multiple panel plots are depicted in Figure~\ref{fig:plot13} and the
single panel plot in \ref{fig:plot2}.

\subsection{Merging and binding}
\label{sec:merge}

As for many rectangular data formats in \proglang{R}, there are
both methods for combining the rows and columns of \code{"zoo"}
objects respectively. For the \code{rbind} method the number of
columns of the combined objects has to be identical and the
indexes may not overlap.
\begin{Schunk}
\begin{Sinput}
R> rbind(z1[5:10], z1[2:3])
\end{Sinput}
\begin{Soutput}
 2004-01-14  2004-01-19  2004-01-27  2004-02-07  2004-02-12  2004-02-16 
 0.02107873 -0.29823529  1.94078850  1.27384445  0.22170438 -2.07607585 
 2004-02-20  2004-02-24 
-1.78439244 -0.19533304 
\end{Soutput}
\end{Schunk}
The \code{c} method simply calls \code{rbind} and hence behaves in the same way.

The \code{cbind} method by default combines the columns by the union of
the indexes and fills the created gaps by \code{NA}s.
\begin{Schunk}
\begin{Sinput}
R> cbind(z1, z2)
\end{Sinput}
\begin{Soutput}
           z1          z2         
2004-01-03          NA  0.94306673
2004-01-05  0.74675994 -0.04149429
2004-01-14  0.02107873          NA
2004-01-17          NA  0.59448077
2004-01-19 -0.29823529 -0.52575918
2004-01-24          NA -0.96739776
2004-01-25  0.68625772          NA
2004-01-27  1.94078850          NA
2004-02-07  1.27384445          NA
2004-02-08          NA  0.95605566
2004-02-12  0.22170438 -0.62733473
2004-02-13          NA -0.92845336
2004-02-16 -2.07607585          NA
2004-02-20 -1.78439244          NA
2004-02-24 -0.19533304          NA
2004-02-25          NA  0.56060280
2004-02-26          NA  0.08291711
\end{Soutput}
\end{Schunk}
In fact, the \code{cbind} method is synonymous with the \code{merge}
method\footnote{Note, that in some situations the column naming in the
resulting object is somewhat problematic in the \code{cbind} method
and the \code{merge} method might provide better formatting of the
column names.}
except that the latter provides additional arguments
which allow for combining the columns by the intersection
of the indexes using the argument \code{all = FALSE}
\begin{Schunk}
\begin{Sinput}
R> merge(z1, z2, all = FALSE)
\end{Sinput}
\begin{Soutput}
           z1          z2         
2004-01-05  0.74675994 -0.04149429
2004-01-19 -0.29823529 -0.52575918
2004-02-12  0.22170438 -0.62733473
\end{Soutput}
\end{Schunk}
Additionally, the filling pattern can be changed in \code{merge},
the naming of the
columns can be modified and the return class of the result can
be specified. In the case of merging of objects with 
different index classes, \proglang{R} gives a warning and tries to
coerce the indexes. Merging objects with different index classes is
generally discouraged---if it is used nevertheless, it is the
responsibility of the user to ensure that the result is as intended.
If at least one of the merged/binded objects was a \code{"zooreg"} 
object, then \code{merge} tries to return a \code{"zooreg"}
object. This is done by assessing whether there is a common maximal 
frequency and by checking whether the resulting index is still
(weakly) regular.

If non-\code{"zoo"} objects are included in merging,
then \code{merge} gives plain vectors/factors/matrices the index of the
first argument (if it is of the same length). Scalars are always added for
the full index without missing values.

\begin{Schunk}
\begin{Sinput}
R> merge(z1, pi, 1:10)
\end{Sinput}
\begin{Soutput}
           z1          pi          1:10       
2004-01-05  0.74675994  3.14159265  1.00000000
2004-01-14  0.02107873  3.14159265  2.00000000
2004-01-19 -0.29823529  3.14159265  3.00000000
2004-01-25  0.68625772  3.14159265  4.00000000
2004-01-27  1.94078850  3.14159265  5.00000000
2004-02-07  1.27384445  3.14159265  6.00000000
2004-02-12  0.22170438  3.14159265  7.00000000
2004-02-16 -2.07607585  3.14159265  8.00000000
2004-02-20 -1.78439244  3.14159265  9.00000000
2004-02-24 -0.19533304  3.14159265 10.00000000
\end{Soutput}
\end{Schunk}

Another function which performs operations along a subset of indexes
is \code{aggregate}, which is discussed in this section although
it does not combine several objects. Using the \code{aggregate} method, \code{"zoo"} objects
are split into subsets along a coarser index grid,
summary statistics are computed for each and then the 
reduced object is returned. In the following example,
first a function is set up which returns for a given \code{"Date"}
value the corresponding first of the month. This function is then
used to compute the coarser grid for the \code{aggregate} call: in
the first example, the grouping is computed explicitely by \verb/firstofmonth(Z.index)/
and the mean of the observations in the month
is returned---in the second example, only the function that computes 
the grouping (when applied to \verb/index(Z)/) is supplied and
the first observation is used for aggregation.

\begin{Schunk}
\begin{Sinput}
R> firstofmonth <- function(x) as.Date(sub("..$", "01", format(x)))
R> aggregate(Z, firstofmonth(Z.index), mean)
\end{Sinput}
\begin{Soutput}
           Aa          Bb          Cc         
2004-02-01  0.53820841  0.04508597 -0.12412352
2004-03-01 -1.18080051  0.58156655 -0.45730045
\end{Soutput}
\begin{Sinput}
R> aggregate(Z, firstofmonth(Z.index), head, 1)
\end{Sinput}
\begin{Soutput}
           Aa         Bb         Cc        
2004-02-01  1.2554339  0.6815732 -0.6329205
2004-03-01 -1.4945833  1.3234122 -1.4944227
\end{Soutput}
\end{Schunk}

\subsection{Mathematical operations}
\label{sec:Ops}

To allow for standard mathematical operations among \code{"zoo"}
objects, \pkg{zoo} extends group generic functions \code{Ops}.
These perform the operations only for the intersection of the
indexes of the objects. As an example, the summation and logical
comparison with $<$ of \code{z1} and \code{z2} yield
\begin{Schunk}
\begin{Sinput}
R> z1 + z2
\end{Sinput}
\begin{Soutput}
2004-01-05 2004-01-19 2004-02-12 
 0.7052657 -0.8239945 -0.4056304 
\end{Soutput}
\begin{Sinput}
R> z1 < z2
\end{Sinput}
\begin{Soutput}
2004-01-05 2004-01-19 2004-02-12 
     FALSE      FALSE      FALSE 
\end{Soutput}
\end{Schunk}

Additionally, methods for transposing \code{t} of \code{"zoo"}
objects---which coerces to a matrix before---and 
computing cumulative quantities such as
\code{cumsum}, \code{cumprod}, \code{cummin}, \code{cummax}
which are all applied column wise.
\begin{Schunk}
\begin{Sinput}
R> cumsum(Z)
\end{Sinput}
\begin{Soutput}
           Aa         Bb         Cc        
2004-02-02  1.2554339  0.6815732 -0.6329205
2004-02-08 -0.2391494  2.0049854 -2.1273432
2004-02-09 -2.1137718  1.1316925 -1.5000035
2004-02-21 -2.2591579  1.5840415 -1.6459775
2004-02-22 -2.0337337  2.1224309 -1.4146162
2004-02-29 -0.8267785  2.4405731 -1.4259082
2004-03-05 -2.0353888  3.8643710 -2.2420530
2004-03-10 -2.1457844  5.2121135 -1.2868283
2004-03-14 -1.3037606  2.4736933 -1.0553214
2004-03-20 -1.4939516  2.5967820 -2.5739429
\end{Soutput}
\end{Schunk}

\subsection{Extracting and replacing the data and the index}
\label{sec:window}

\pkg{zoo} provides several generic functions and methods
to work on the data contained in a \code{"zoo"} object, the
index (or time) attribute associated to it, and on both data and
index.

The data stored in \code{"zoo"} objects can be extracted by
\code{coredata} which strips off all \code{"zoo"}-specific attributes and 
it can be replaced using \code{coredata<-}. Both are new generic
functions\footnote{The \code{coredata} functionality is similar in spirit to the \code{core}
function in \pkg{its} and \code{value} in \pkg{tseries}. However, the 
focus of those functions is somewhat narrower and we try to provide 
more general purpose generic functions. See the respective manual
page for more details.}
with methods for \code{"zoo"} objects as illustrated in the following
example.
\begin{Schunk}
\begin{Sinput}
R> coredata(z1)
\end{Sinput}
\begin{Soutput}
 [1]  0.74675994  0.02107873 -0.29823529  0.68625772  1.94078850  1.27384445
 [7]  0.22170438 -2.07607585 -1.78439244 -0.19533304
\end{Soutput}
\begin{Sinput}
R> coredata(z1) <- 1:10
R> z1
\end{Sinput}
\begin{Soutput}
2004-01-05 2004-01-14 2004-01-19 2004-01-25 2004-01-27 2004-02-07 2004-02-12 
         1          2          3          4          5          6          7 
2004-02-16 2004-02-20 2004-02-24 
         8          9         10 
\end{Soutput}
\end{Schunk}

The index associated with a \code{"zoo"} object can be extracted
by \code{index} and modified by \code{index<-}. As the interpretation
of the index as ``time'' in time series applications is natural,
there are also synonymous methods \code{time} and \code{time<-}. 
Hence, the commands \code{index(z2)} and \code{time(z2)}
return equivalent results.
\begin{Schunk}
\begin{Sinput}
R> index(z2)
\end{Sinput}
\begin{Soutput}
 [1] "2004-01-03 Eastern Standard Time" "2004-01-05 Eastern Standard Time"
 [3] "2004-01-17 Eastern Standard Time" "2004-01-19 Eastern Standard Time"
 [5] "2004-01-24 Eastern Standard Time" "2004-02-08 Eastern Standard Time"
 [7] "2004-02-12 Eastern Standard Time" "2004-02-13 Eastern Standard Time"
 [9] "2004-02-25 Eastern Standard Time" "2004-02-26 Eastern Standard Time"
\end{Soutput}
\end{Schunk}
The index scale of \code{z2} can be changed to that of \code{z1} by
\begin{Schunk}
\begin{Sinput}
R> index(z2) <- index(z1)
R> z2
\end{Sinput}
\begin{Soutput}
 2004-01-05  2004-01-14  2004-01-19  2004-01-25  2004-01-27  2004-02-07 
 0.94306673 -0.04149429  0.59448077 -0.52575918 -0.96739776  0.95605566 
 2004-02-12  2004-02-16  2004-02-20  2004-02-24 
-0.62733473 -0.92845336  0.56060280  0.08291711 
\end{Soutput}
\end{Schunk}

The start and the end of the index/time vector can be queried by
\code{start} and \code{end}:
\begin{Schunk}
\begin{Sinput}
R> start(z1)
\end{Sinput}
\begin{Soutput}
[1] "2004-01-05 Eastern Standard Time"
\end{Soutput}
\begin{Sinput}
R> end(z1)
\end{Sinput}
\begin{Soutput}
[1] "2004-02-24 Eastern Standard Time"
\end{Soutput}
\end{Schunk}

To work on both data and index/time, \pkg{zoo} provides
\code{window} and \code{window<-} methods for \code{"zoo"} objects.
In both cases the window is specified by
\begin{Scode}
window(x, index, start, end)
\end{Scode}
where \code{x} is the \code{"zoo"} object, \code{index} is a set
of indexes to be selected (by default the full index of \code{x})
and \code{start} and \code{end} can be used to restrict the 
\code{index} set. 
\begin{Schunk}
\begin{Sinput}
R> window(Z, start = as.Date("2004-03-01"))
\end{Sinput}
\begin{Soutput}
           Aa         Bb         Cc        
2004-03-05 -1.2086102  1.4237978 -0.8161448
2004-03-10 -0.1103956  1.3477425  0.9552247
2004-03-14  0.8420238 -2.7384202  0.2315069
2004-03-20 -0.1901910  0.1230887 -1.5186216
\end{Soutput}
\begin{Sinput}
R> window(Z, index = index(Z)[5:8], end = as.Date("2004-03-01"))
\end{Sinput}
\begin{Soutput}
           Aa          Bb          Cc         
2004-02-22  0.22542418  0.53838938  0.23136133
2004-02-29  1.20695518  0.31814222 -0.01129202
\end{Soutput}
\end{Schunk}

The first example selects all observations starting from 2004-03-01
whereas the second selects from the from the 5th to 8th observation
those up to 2004-03-01.

The same syntax can be used for the corresponding replacement function.
\begin{Schunk}
\begin{Sinput}
R> window(z1, end = as.POSIXct("2004-02-01")) <- 9:5
R> z1
\end{Sinput}
\begin{Soutput}
2004-01-05 2004-01-14 2004-01-19 2004-01-25 2004-01-27 2004-02-07 2004-02-12 
         9          8          7          6          5          6          7 
2004-02-16 2004-02-20 2004-02-24 
         8          9         10 
\end{Soutput}
\end{Schunk}

Two methods that are standard in time series applications
are \code{lag} and \code{diff}. These are available with the same
arguments as the \code{"ts"} methods.\footnote{\code{diff} also
has an additional argument that also allows for geometric and
not only allows arithmetic differences. Furthermore, note the sign
of the lag in \code{lag}: by default it is positive and shifts the 
observations \emph{forward}, to obtain the more standard \emph{backward}
shift the lag has to be negative.}

\begin{Schunk}
\begin{Sinput}
R> lag(z1, k = -1)
\end{Sinput}
\begin{Soutput}
2004-01-14 2004-01-19 2004-01-25 2004-01-27 2004-02-07 2004-02-12 2004-02-16 
         9          8          7          6          5          6          7 
2004-02-20 2004-02-24 
         8          9 
\end{Soutput}
\begin{Sinput}
R> merge(z1, lag(z1, k = 1))
\end{Sinput}
\begin{Soutput}
           z1 lag(z1, k = 1)
2004-01-05  9  8            
2004-01-14  8  7            
2004-01-19  7  6            
2004-01-25  6  5            
2004-01-27  5  6            
2004-02-07  6  7            
2004-02-12  7  8            
2004-02-16  8  9            
2004-02-20  9 10            
2004-02-24 10 NA            
\end{Soutput}
\begin{Sinput}
R> diff(z1)
\end{Sinput}
\begin{Soutput}
2004-01-14 2004-01-19 2004-01-25 2004-01-27 2004-02-07 2004-02-12 2004-02-16 
        -1         -1         -1         -1          1          1          1 
2004-02-20 2004-02-24 
         1          1 
\end{Soutput}
\end{Schunk}

\subsection[Coercion to and from "zoo"]{Coercion to and from \code{"zoo"}}
\label{sec:as.zoo}

Coercion to and from \code{"zoo"} objects is available for objects of
various classes, in particular \code{"ts"}, \code{"irts"} and \code{"its"}
objects can be coerced to \code{"zoo"} and back if the index is of the appropriate
class.\footnote{Coercion from \code{"zoo"} to \code{"irts"} is contained in the
\pkg{tseries} package.}

Coercion between \code{"zooreg"} and \code{"zoo"} is also available and is essentially
dropping the \code{"frequency"} attribute or trying to add one, respectively.

Furthermore, \code{"zoo"} objects can be coerced to vectors, matrices, lists and
data frames (the latter dropping the index/time attribute). A simple example is
\begin{Schunk}
\begin{Sinput}
R> as.data.frame(Z)
\end{Sinput}
\begin{Soutput}
                   Aa         Bb          Cc
2004-02-02  1.2554339  0.6815732 -0.63292049
2004-02-08 -1.4945833  1.3234122 -1.49442269
2004-02-09 -1.8746225 -0.8732929  0.62733971
2004-02-21 -0.1453861  0.4523490 -0.14597401
2004-02-22  0.2254242  0.5383894  0.23136133
2004-02-29  1.2069552  0.3181422 -0.01129202
2004-03-05 -1.2086102  1.4237978 -0.81614483
2004-03-10 -0.1103956  1.3477425  0.95522468
2004-03-14  0.8420238 -2.7384202  0.23150695
2004-03-20 -0.1901910  0.1230887 -1.51862157
\end{Soutput}
\end{Schunk}

\subsection[NA handling]{\code{NA} handling}
\label{sec:NA}

Four methods for dealing with \code{NA}s (missing observations) 
in the observations are applicable to \code{"zoo"} objects:
\code{na.omit}, \code{na.contiguous}, \code{na.approx} and \code{na.locf}.
\code{na.omit}---or its default method to be more precise---returns a \code{"zoo"}
object with incomplete observations removed. \code{na.contiguous}
extracts the longest consecutive stretch of non-missing values.
Furthermore, new generic functions
\code{na.approx} and \code{na.locf} and corresponding default methods are introduced in \pkg{zoo}.
The former replaces \code{NA}s by linear interpolation (using the
function \code{approx}) and the name of the latter
stands for \underline{l}ast \underline{o}bservation \underline{c}arried
\underline{f}orward. It replaces missing observations by the most recent
non-\code{NA} prior to it. Leading \code{NA}s, which cannot be replaced
by previous observations, are removed in both functions by default.

\begin{Schunk}
\begin{Sinput}
R> z1[sample(1:10, 3)] <- NA
R> z1
\end{Sinput}
\begin{Soutput}
2004-01-05 2004-01-14 2004-01-19 2004-01-25 2004-01-27 2004-02-07 2004-02-12 
         9         NA          7          6          5          6         NA 
2004-02-16 2004-02-20 2004-02-24 
         8          9         NA 
\end{Soutput}
\begin{Sinput}
R> na.omit(z1)
\end{Sinput}
\begin{Soutput}
2004-01-05 2004-01-19 2004-01-25 2004-01-27 2004-02-07 2004-02-16 2004-02-20 
         9          7          6          5          6          8          9 
\end{Soutput}
\begin{Sinput}
R> na.contiguous(z1)
\end{Sinput}
\begin{Soutput}
2004-01-19 2004-01-25 2004-01-27 2004-02-07 
         7          6          5          6 
\end{Soutput}
\begin{Sinput}
R> na.approx(z1)
\end{Sinput}
\begin{Soutput}
2004-01-05 2004-01-14 2004-01-19 2004-01-25 2004-01-27 2004-02-07 2004-02-12 
  9.000000   7.714286   7.000000   6.000000   5.000000   6.000000   7.111111 
2004-02-16 2004-02-20 
  8.000000   9.000000 
\end{Soutput}
\begin{Sinput}
R> na.approx(z1, 1:NROW(z1))
\end{Sinput}
\begin{Soutput}
2004-01-05 2004-01-14 2004-01-19 2004-01-25 2004-01-27 2004-02-07 2004-02-12 
         9          8          7          6          5          6          7 
2004-02-16 2004-02-20 
         8          9 
\end{Soutput}
\begin{Sinput}
R> na.locf(z1)
\end{Sinput}
\begin{Soutput}
2004-01-05 2004-01-14 2004-01-19 2004-01-25 2004-01-27 2004-02-07 2004-02-12 
         9          9          7          6          5          6          6 
2004-02-16 2004-02-20 2004-02-24 
         8          9          9 
\end{Soutput}
\end{Schunk}

As the above example illustrates, \code{na.approx} uses by default
the underlying time scale for interpolation. This can be changed, e.g.,
to an equidistant spacing, by setting the second argument of
\code{na.approx}.

\subsection{Rolling functions}
\label{sec:rolling}

A typical task to be performed on ordered observations is to evaluate some
function, e.g., computing the mean, in a window of observations that is moved
over the full sample period. The resulting statistics are usually synonymously referred to
as rolling/running/moving statistics. In \pkg{zoo}, the generic function \code{rapply} 
is provided along with a \code{"zoo"} and a \code{"ts"} method. The most important arguments
are

\begin{Scode}
rapply(data, width, FUN)
\end{Scode}

where the function \code{FUN} is applied to a rolling window of size \code{width}
of the observations \code{data}. The function \code{rapply} currently only evaluates
the function for windows of full size \code{width}, hence the result has \code{width - 1}
fewer observations than the original series. But it can be determined whether the `lost'
observations should be padded with \code{NA}s and whether the result should be left-
or right-aligned or centered (default) with respect to the original index.

\begin{Schunk}
\begin{Sinput}
R> rapply(Z, 5, sd)
\end{Sinput}
\begin{Soutput}
           Aa        Bb        Cc       
2004-02-09 1.2814876 0.8018950 0.8218959
2004-02-21 1.2658555 0.7891358 0.8025043
2004-02-22 1.2102011 0.8206819 0.5319727
2004-02-29 0.8662296 0.5266261 0.6411751
2004-03-05 0.9363400 1.7011273 0.6356144
2004-03-10 0.9508642 1.6892246 0.9578196
\end{Soutput}
\begin{Sinput}
R> rapply(Z, 5, sd, na.pad = TRUE, align = "left")
\end{Sinput}
\begin{Soutput}
           Aa        Bb        Cc       
2004-02-02 1.2814876 0.8018950 0.8218959
2004-02-08 1.2658555 0.7891358 0.8025043
2004-02-09 1.2102011 0.8206819 0.5319727
2004-02-21 0.8662296 0.5266261 0.6411751
2004-02-22 0.9363400 1.7011273 0.6356144
2004-02-29 0.9508642 1.6892246 0.9578196
2004-03-05        NA        NA        NA
2004-03-10        NA        NA        NA
2004-03-14        NA        NA        NA
2004-03-20        NA        NA        NA
\end{Soutput}
\end{Schunk}

To improve the performance of \code{rapply(x, k, }\textit{foo}\code{)} for some frequently
used functions \textit{foo}, more efficient implementations \code{roll}\textit{foo}\code{(x, k)}
are available (and also called by \code{rapply}). 
Currently, these are the generic functions \code{rollmean}, \code{rollmedian}
and \code{rollmax} which have methods for \code{"zoo"} and \code{"ts"} series and a 
default method for plain vectors.

\begin{Schunk}
\begin{Sinput}
R> rollmean(z2, 5, na.pad = TRUE)
\end{Sinput}
\begin{Soutput}
   2004-01-05    2004-01-14    2004-01-19    2004-01-25    2004-01-27 
           NA            NA  0.0005792538  0.0031770388 -0.1139910497 
   2004-02-07    2004-02-12    2004-02-16    2004-02-20    2004-02-24 
-0.4185778750 -0.2013054791  0.0087574946            NA            NA 
\end{Soutput}
\end{Schunk}

\section[Combining zoo with other packages]{Combining \pkg{zoo} with other packages}
\label{sec:combining}

The main purpose of the package \pkg{zoo} is to provide basic infrastructure for
working with indexed totally ordered observations that can be either employed by
users directly or can be a basic ingredient on top of which other packages can
build. The latter is illustrated with a few brief examples involving the packages
\pkg{strucchange}, \pkg{tseries} and \pkg{fCalendar} in this section. Finally, the 
classes \code{"yearmon"} and \code{"yearqtr"} (provided in \pkg{zoo})
are used for illustrating how \pkg{zoo} can be extended by creating a new index class.

\subsection[strucchange: Empirical fluctuation processes]{\pkg{strucchange}: Empirical fluctuation processes}
\label{sec:strucchange}

The package \pkg{strucchange} provides a collection of methods for testing,
monitoring and dating structural changes, in particular in linear regression models.
Tests for structural change assess whether the parameters of a model remain
constant over an ordering with respect to a specified variable, usually time.
To adequatly store and visualize empirical fluctuation processes which 
capture instabilities over this ordering, a data type for indexed ordered
observations is required. This was the motivation for starting the \pkg{zoo}
project.

A simple example for the need of \code{"zoo"} objects in \pkg{strucchange}
which can not be (easily) implemented by other irregular time series classes
available in \proglang{R} is described in the following. We assess the constancy of the
electrical resistance over the apparent juice content of kiwi fruits.\footnote{A different
approach would be to test whether the slope of a regression of electrical resistance
on juice content changes with increasing juice content, i.e., to test for
instabilities in \code{ohms \~{} juice} instead of \code{ohms \~{} 1}. Both lead to 
similar results.} The data
set \code{fruitohms} is contained in the \pkg{DAAG} package \citep{zoo:DAAG:2004}.
The fitted \code{ocus} object contains the OLS-based CUSUM process for the mean
of the electrical resistance (variable \code{ohms}) indexed by the juice
content (variable \code{juice}).

\begin{Schunk}
\begin{Sinput}
R> library(strucchange)
R> library(DAAG)
R> data(fruitohms)
R> ocus <- gefp(ohms ~ 1, order.by = ~juice, data = fruitohms)
\end{Sinput}
\end{Schunk}

\begin{figure}[h!]
\begin{center}
\begin{Schunk}
\begin{Sinput}
R> plot(ocus)
\end{Sinput}
\end{Schunk}
\includegraphics{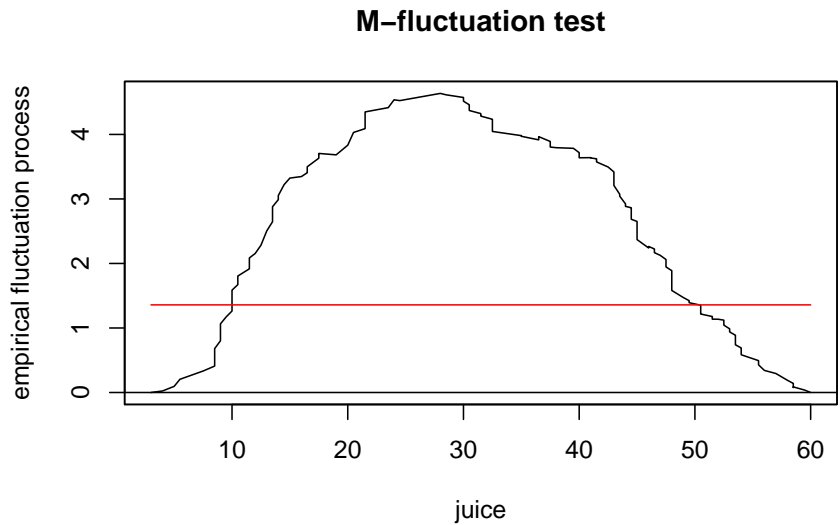}
\caption{\label{fig:strucchange} Empirical M-fluctuation process for \code{fruitohms} data}
\end{center}
\end{figure}

This OLS-based CUSUM process can be visualized using the \code{plot} method
for \code{"gefp"} objects which builds on the \code{"zoo"} method and yields in
this case the plot in Figure~\ref{fig:strucchange} showing the process which
crosses its 5\% critical value and 
thus signals a significant decrease in the mean electrical resistance over the
juice content. For more information on the package \pkg{strucchange} and the 
function \code{gefp} see \cite{zoo:Zeileis+Leisch+Hornik:2002} and 
\cite{zoo:Zeileis:2004}.

\subsection[tseries: Historical financial data]{\pkg{tseries}: Historical financial data}
\label{sec:tseries}

A typical application for irregular time series which became increasingly
important over the last years in computational statistics and finance is
daily (or higher frequent) financial data. The package \pkg{tseries} provides
the function \code{get.hist.quote} for obtaining historical financial data
by querying Yahoo! Finance at \url{http://finance.yahoo.com/},
an online portal quoting data provided by Reuters. The following code
queries the quotes of Lucent Technologies starting from 2001-01-01
until 2004-09-30:

\begin{Schunk}
\begin{Sinput}
R> library(tseries)
R> LU <- get.hist.quote(instrument = "LU", start = "2001-01-01", 
+     end = "2004-09-30", origin = "1970-01-01")
\end{Sinput}
\end{Schunk}

\begin{figure}[h!]
\begin{center}
\begin{Schunk}
\begin{Sinput}
R> plot(diff(log(LU)))
\end{Sinput}
\end{Schunk}
\includegraphics{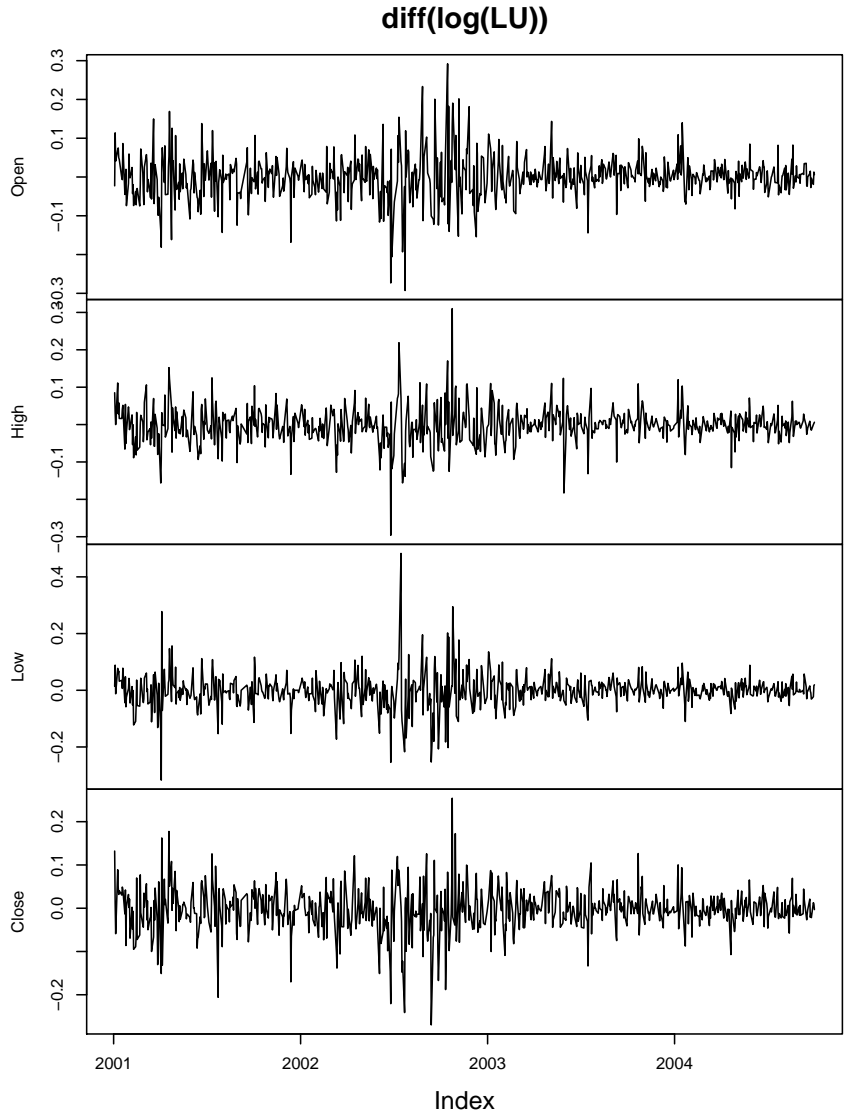}
\caption{\label{fig:tseries} Log-difference returns for Lucent Technologies}
\end{center}
\end{figure}

In the returned \code{LU} object the irregular data is stored by extending
it in a regular grid and filling the gaps with \code{NA}s. The time is stored
in days starting from an \code{origin}, in this case specified to be 1970-01-01, the
origin used by the \code{Date} class.
This series can be transformed easily into an irregular \code{"zoo"} series 
using a \code{"Date"} index. The log-difference returns for Lucent 
Technologies is depicted in Figure~\ref{fig:tseries}.

\begin{Schunk}
\begin{Sinput}
R> LU <- as.zoo(LU)
R> index(LU) <- as.Date(index(LU))
R> LU <- na.omit(LU)
\end{Sinput}
\end{Schunk}

\subsection[fCalendar: Indexes of class "timeDate"]{\pkg{fCalendar}: Indexes of class \code{"timeDate"}}
\label{sec:fCalendar}

Although the methods in \pkg{zoo} work out of the box for many index classes,
it might be necessary for some index classes to provide \code{c}, \code{length},
\code{ORDER} and \code{MATCH} methods such that the methods in \pkg{zoo} 
work properly. An example for such an index class which requires a bit more
attention is \code{"timeDate"} from the \pkg{fCalendar} package.

But after the necessary methods have been defined
\begin{Schunk}
\begin{Sinput}
R> length.timeDate <- function(x) prod(x@Dim)
R> ORDER.timeDate <- function(x, ...) order(as.POSIXct(x), ...)
R> MATCH.timeDate <- function(x, table, nomatch = NA, ...) match(as.POSIXct(x), 
+     as.POSIXct(table), nomatch = NA, ...)
\end{Sinput}
\end{Schunk}
the class \code{"timeDate"} can be used for indexing \code{"zoo"} objects.
The following example illustrates how \code{z2} can be transformed
to use the \code{"timeDate"} class.
\begin{Schunk}
\begin{Sinput}
R> library(fCalendar)
R> z2td <- zoo(coredata(z2), timeDate(index(z2), FinCenter = "GMT"))
R> z2td
\end{Sinput}
\begin{Soutput}
 2004-01-05  2004-01-14  2004-01-19  2004-01-25  2004-01-27  2004-02-07 
 0.94306673 -0.04149429  0.59448077 -0.52575918 -0.96739776  0.95605566 
 2004-02-12  2004-02-16  2004-02-20  2004-02-24 
-0.62733473 -0.92845336  0.56060280  0.08291711 
\end{Soutput}
\end{Schunk}

\subsection[The classes "yearmon" and "yearqtr": Roll your own index]{The classes \code{"yearmon"} and \code{"yearqtr"}: Roll your own index}
\label{sec:yearmon}

One of the strengths of the \pkg{zoo} package is its independence of the
index class, such that the index can be easily customized. The previous section
already explained how an existing class (\code{"timeDate"}) can be used as
the index if the necessary methods are created. This section has a similar
but slightly different focus: it describes how new index classes can be created
addressing a certain type of indexes. These classes are \code{"yearmon"} and
\code{"yearqtr"} (already contained in \pkg{zoo}) which provide indexes for
monthly and quarterly data respectively.
As the code is virtually identical for both classes---except that one has the 
frequency 12 and the other 4---we will only discuss \code{"yearmon"} explicitly.

Of course, monthly data can simply be stored using a numeric index just
as the class \code{"ts"} does. The problem is that this does not have the meta-information
attached that this is really specifying monthly data which is in \code{"yearmon"}
simply added by a class attribute. Hence, the class creator is simply defined as

\begin{Scode}
yearmon <- function(x) structure(floor(12*x + .0001)/12, class = "yearmon")
\end{Scode}

which is very similar to the \code{as.yearmon} coercion functions provided.

As \code{"yearmon"} data is now explicitly declared to describe monthly data,
this can be exploited for coercion to other time classes: either to coarser time scales
such as \code{"yearqtr"} or to finer time scales such as
\code{"Date"}, \code{"POSIXct"} or \code{"POSIXlt"} which by default associate the first day
within a month with a \code{"yearmon"} observation. Adding a \code{format}
and \code{as.character} method produces human readable character representations
of \code{"yearmon"} data and \code{Ops} and \code{MATCH} methods complete the
methods needed for conveniently  working with monthly data in \pkg{zoo}. Note,
that all of these methods are  very simple and rather obvious (as can be seen in
the \pkg{zoo} sources), but prove very helpful in the following examples.

First, we create a regular series \code{zr3} with \code{"yearmon"} index which
leads to improved printing compared to the regular series \code{zr1} and \code{zr2}
from Section~\ref{sec:zooreg}.

\begin{Schunk}
\begin{Sinput}
R> zr3 <- zooreg(rnorm(9), start = yearmon(2000), frequency = 12)
R> zr3
\end{Sinput}
\begin{Soutput}
   Jan 2000    Feb 2000    Mar 2000    Apr 2000    May 2000    Jun 2000 
-0.30969096  0.08699142 -0.64837101 -0.62786277 -0.61932674 -0.95506154 
   Jul 2000    Aug 2000    Sep 2000 
-1.91736406  0.38108885  1.51405511 
\end{Soutput}
\end{Schunk}

This could be aggregated to quarterly data via

\begin{Schunk}
\begin{Sinput}
R> aggregate(zr3, as.yearqtr, mean)
\end{Sinput}
\begin{Soutput}
   2000 Q1    2000 Q2    2000 Q3 
-0.2903569 -0.7340837 -0.0074067 
\end{Soutput}
\end{Schunk}

The index can easily be transformed to \code{"Date"}, the default being the first day
of the month but which can also be changed to the last day of the month.

\begin{Schunk}
\begin{Sinput}
R> as.Date(index(zr3))
\end{Sinput}
\begin{Soutput}
[1] "2000-01-01" "2000-02-01" "2000-03-01" "2000-04-01" "2000-05-01"
[6] "2000-06-01" "2000-07-01" "2000-08-01" "2000-09-01"
\end{Soutput}
\begin{Sinput}
R> as.Date(index(zr3), frac = 1)
\end{Sinput}
\begin{Soutput}
[1] "2000-01-31" "2000-02-29" "2000-03-31" "2000-04-30" "2000-05-31"
[6] "2000-06-30" "2000-07-31" "2000-08-31" "2000-09-30"
\end{Soutput}
\end{Schunk}

Furthermore, \code{"yearmon"} indexes can easily be coerced to \code{"POSIXct"} such
that the series could be exported as a \code{"its"} or \code{"irts"} series.

\begin{Schunk}
\begin{Sinput}
R> index(zr3) <- as.POSIXct(index(zr3))
R> as.irts(zr3)
\end{Sinput}
\begin{Soutput}
2000-01-01 00:00:00 GMT -0.3097
2000-02-01 00:00:00 GMT 0.08699
2000-03-01 00:00:00 GMT -0.6484
2000-04-01 00:00:00 GMT -0.6279
2000-05-01 00:00:00 GMT -0.6193
2000-06-01 00:00:00 GMT -0.9551
2000-07-01 00:00:00 GMT -1.917
2000-08-01 00:00:00 GMT 0.3811
2000-09-01 00:00:00 GMT 1.514
\end{Soutput}
\end{Schunk}

Again, this functionality makes switching between different time scales or index 
representations particularly easy and \pkg{zoo} provides the user with the flexibility
to adjust a certain index to his/her problem of interest.

\newpage

\section{Summary and outlook} \label{sec:summary}

The package \pkg{zoo} provides an \proglang{S3} class and methods
for indexed totally ordered observations, such as both regular and irregular time series.
Its key design goals are independence of a particular index class 
and compatibility with standard generics similar to the behaviour of 
the corresponding \code{"ts"} methods. This paper describes how
these are implemented in \pkg{zoo} and illustrates the usage of 
the methods for plotting, merging and
binding, several mathematical operations, extracting and replacing data
and index, coercion and \code{NA} handling.

An indexed object of class \code{"zoo"} can be thought of as data plus index
where the data are essentially vectors or matrices and the index can be
a vector of (in principle) arbitrary class. For (weakly) regular \code{"zooreg"}
series, a \code{"frequency"} attribute is stored in addition. Therefore, objects of classes
\code{"ts"}, \code{"its"}, \code{"irts"} and \code{"timeSeries"} can easily
be transformed into \code{"zoo"} objects---the reverse transformation is also possible 
provided that the index fulfills the restrictions of the respective class.
Hence, the \code{"zoo"} class can also be used as the basis for other
classes of indexed observations and more specific functionality can be built on
top of it. Furthermore, it bridges the gap between irregular and regular series,
facilitating operations such as \code{NA} handling compared to \code{"ts"}.

Whereas a lot of effort was put into achieving independence of a particular
index class, the types of data that can be indexed with \code{"zoo"} are currently
limited to vectors and matrices, typically containing numeric values. Although,
there is some limited support available for indexed factors, one important 
direction for future development of \pkg{zoo} is to add better support for other
objects that can also naturally be indexed including specifically factors, data
frames and lists.

\section*{Computational details}

The results in this paper were obtained using \proglang{R}
2.1.0 with the packages
\pkg{zoo} 1.0--0,
\pkg{strucchange} 1.2--10,
\pkg{fCalendar} 201.10060,
\pkg{tseries} 0.9--27 and
\pkg{DAAG} 0.46.
\proglang{R} itself and all packages used are available from
CRAN at \url{http://CRAN.R-project.org/}.

\bibliography{zoo}

\begin{thebibliography}{7}
\providecommand{\natexlab}[1]{#1}
\providecommand{\url}[1]{\texttt{#1}}
\expandafter\ifx\csname urlstyle\endcsname\relax
  \providecommand{\doi}[1]{doi: #1}\else
  \providecommand{\doi}{doi: \begingroup \urlstyle{rm}\Url}\fi

\bibitem[Heywood(2004)]{zoo:its:2004}
Giles Heywood.
\newblock \emph{\pkg{its}: Irregular Time Series}.
\newblock Portfolio \& Risk Advisory Group and Commerzbank Securities, 2004.
\newblock \proglang{R} package version 1.0.4.

\bibitem[Maindonald and Braun(2004)]{zoo:DAAG:2004}
John Maindonald and W.~John Braun.
\newblock \emph{\pkg{DAAG}: Data Analysis and Graphics}, 2004.
\newblock URL \url{http://www.stats.uwo.ca/DAAG/}.
\newblock \proglang{R} package version 0.46.

\bibitem[{\proglang{R} Development Core Team}(2005)]{zoo:R:2005}
{\proglang{R} Development Core Team}.
\newblock \emph{\proglang{R}: {A} Language and Environment for Statistical
  Computing}.
\newblock \proglang{R} Foundation for Statistical Computing, Vienna, Austria,
  2005.
\newblock URL \url{http://www.R-project.org/}.
\newblock {ISBN} 3-900051-00-3.

\bibitem[Trapletti(2005)]{zoo:tseries:2005}
Adrian Trapletti.
\newblock \emph{\pkg{tseries}: Time Series Analysis and Computational Finance},
  2005.
\newblock \proglang{R} package version 0.9-25.

\bibitem[Wuertz(2005)]{zoo:fCalendar:2004}
Diethelm Wuertz.
\newblock \emph{\pkg{Rmetrics}: {A}n Environment and Software Collection for
  Teaching Financial Engineering and Computational Finance}, 2005.
\newblock URL \url{http://www.Rmetrics.org/}.
\newblock \proglang{R} package \pkg{fCalendar}, version 201.10059.

\bibitem[Zeileis(2004)]{zoo:Zeileis:2004}
Achim Zeileis.
\newblock Implementing a class of structural change tests: {A}n econometric
  computing approach.
\newblock Report~7, Department of Statistics and Mathematics,
  Wirtschaftsuniversit\"at Wien, Research Report Series, July 2004.
\newblock URL \url{http://epub.wu-wien.ac.at/}.

\bibitem[Zeileis et~al.(2002)Zeileis, Leisch, Hornik, and
  Kleiber]{zoo:Zeileis+Leisch+Hornik:2002}
Achim Zeileis, Friedrich Leisch, Kurt Hornik, and Christian Kleiber.
\newblock \pkg{strucchange}: {A}n \proglang{R} package for testing for
  structural change in linear regression models.
\newblock \emph{Journal of Statistical Software}, 7\penalty0 (2):\penalty0
  1--38, 2002.
\newblock URL \url{http://www.jstatsoft.org/v07/i02/}.

\end{thebibliography}

\newpage

\begin{appendix}
\section{Reference card}
\begin{tabular}{rp{9cm}}
\multicolumn{2}{l}{\textbf{Creation}} \\
\code{zoo(x, order.by)} & creation of a \code{"zoo"} object
  from the observations \code{x} (a vector or a matrix) and an index
  \code{order.by} by which the observations are ordered. \\
& For computations on arbitrary index classes, methods to the 
  following genric functions are assumed to work: combining \code{c()},
  querying length \code{length()}, subsetting \code{[}, ordering
  \code{ORDER()} and value matching \code{MATCH()}. For pretty
  printing an \code{as.character} and/or \code{index2char} method
  might be helpful.\\[0.5cm]

\multicolumn{2}{l}{\textbf{Creation of regular series}} \\
\code{zoo(x, order.by, freq)} & works as above but creates a \code{"zooreg"}
  object which inherits from \code{"zoo"} if the frequency \code{freq} complies
  with the index \code{order.by}. An \code{as.numeric} method has to be
  available for the index class.\\
\code{zooreg(x, start, end, freq)} & creates a \code{"zooreg"} series
  with a numeric index as above and has (almost) the same interface as
  \code{ts()}.\\[0.5cm]

\multicolumn{2}{l}{\textbf{Standard methods}} \\
\code{plot} & plotting \\
\code{lines} & adding a \code{"zoo"} series to a plot \\
\code{print} & printing \\
\code{summary} & summarizing (column-wise) \\
\code{str} & displaying structure of \code{"zoo"} objects \\
\code{head}, \code{tail} & head and tail of \code{"zoo"} objects \\[0.5cm]

\multicolumn{2}{l}{\textbf{Coercion}} \\
\code{as.zoo} & coercion to \code{"zoo"} is available for objects
    of class \code{"ts"}, \code{"its"}, \code{"irts"} (plus a default
    method).\\
\code{as.}\textit{class}\code{.zoo} & coercion from \code{"zoo"} to
    other classes. Currently available for \textit{class} in \code{"matrix"},
    \code{"vector"}, \code{"data.frame"}, \code{"list"}, \code{"irts"},
    \code{"its"} and \code{"ts"}. \\
\code{is.zoo} & querying wether an object is of class \code{"zoo"} \\[0.5cm]

\multicolumn{2}{l}{\textbf{Merging and binding}} \\
\code{merge} & union, intersection, left join, right join along indexes\\
\code{cbind} & column binding along the intersection of the index\\
\code{c}, \code{rbind} & combining/row binding (indexes may not overlap)\\
\code{aggregate} & compute summary statistics along a coarser grid of indexes \\[0.5cm]

\multicolumn{2}{l}{\textbf{Mathematical operations}} \\
\code{Ops} & group generic functions performed along the intersection of indexes\\
\code{t} & transposing (coerces to \code{"matrix"} before) \\
\code{cumsum} & compute (columnwise) cumulative quantities: sums
    \code{cumsum()}, products \code{cumprod()}, maximum \code{cummax()},
    minimum \code{cummin()}.\\[0.5cm]
\end{tabular}

\newpage

\begin{tabular}{rp{9cm}}
\multicolumn{2}{l}{\textbf{Extracting and replacing data and index}} \\
\code{index, time} & extract the index of a series\\
\code{index<-}, \code{time<-} & replace the index of a series\\
\code{coredata}, \code{coredata<-} & extract and replace the data associated with a \code{"zoo"} object\\
\code{lag} & lagged observations \\
\code{diff} & arithmetic and geometric differences \\
\code{start, end} & querying start and end of a series \\
\code{window, window<-} & subsetting of \code{"zoo"} objects
    using their index\\[0.5cm]

\multicolumn{2}{l}{\textbf{\code{NA} handling}} \\
\code{na.omit} & omit \code{NA}s \\
\code{na.contiguous} & compute longest sequence of non-\code{NA} observations \\
\code{na.locf} & impute \code{NA}s by carrying forward the last observation\\
\code{na.approx} & impute \code{NA}s by interpolation\\[0.5cm]

\multicolumn{2}{l}{\textbf{Rolling functions}} \\
\code{rapply} & apply a function to rolling margin of an array \\
\code{rollmean} & more efficient functions for computing the rolling mean, median
  and maximum are \code{rollmean()}, \code{rollmedian()} and \code{rollmax()}, respectively\\[0.5cm]

\multicolumn{2}{l}{\textbf{Methods for regular series}} \\
\code{is.regular} & checks whether a series is weakly (or strictly if \code{strict = TRUE})
  regular \\
\code{frequency}, \code{deltat} & extracts the frequency or its reciprocal value
  respectively from a series, for \code{"zoo"} series the functions try to determine
  the regularity and frequency in a data-driven way\\
\code{cycle} & gives the position in the cycle of a regular series \\[0.5cm]

\end{tabular}

\end{appendix}

\end{document}